# STRONG MEMORYLESS TIMES AND RARE EVENTS IN MARKOV RENEWAL POINT PROCESSES

By Torkel Erhardsson

*Royal Institute of Technology*

Let $W$ be the number of points in $(0, t]$ of a stationary finite-state Markov renewal point process. We derive a bound for the total variation distance between the distribution of $W$ and a compound Poisson distribution. For any nonnegative random variable $\zeta$, we construct a "strong memoryless time" $\hat{\zeta}$ such that $\zeta - t$ is exponentially distributed conditional on $\{\hat{\zeta} \leq t, \zeta > t\}$, for each $t$. This is used to embed the Markov renewal point process into another such process whose state space contains a frequently observed state which represents loss of memory in the original process. We then write $W$ as the accumulated reward of an embedded renewal reward process, and use a compound Poisson approximation error bound for this quantity by Erhardsson. For a renewal process, the bound depends in a simple way on the first two moments of the interrenewal time distribution, and on two constants obtained from the Radon–Nikodym derivative of the interrenewal time distribution with respect to an exponential distribution. For a Poisson process, the bound is 0.

**1. Introduction.** In this paper, we are concerned with rare events in stationary finite-state Markov renewal point processes (MRPPs). An MRPP is a marked point process on $\mathbb{R}$ or $\mathbb{Z}$ (continuous or discrete time). Each point of an MRPP has an associated mark, or state. The distance in time between two successive points and the state of the second point are jointly conditionally independent of the past given the state of the first point. A renewal process is a special case of an MRPP, and any finite-state Markov or semi-Markov process can be constructed using a suitable MRPP, simply by defining the state of the process at time $t$ to be the state of the most recently observed point of the MRPP.

The number of points of a stationary MRPP in $(0, t]$ with states in a certain subset $B$ of the state space is an important quantity in many applications. For example, the number of visits to $B$ in $(0, t]$ by a stationary Markov chain can be expressed in this way. If points with states in $B$ are rare, this quantity should be approximately compound Poisson distributed. Heuristically, the set of such points can be partitioned into disjoint clumps, the sizes of which are approximately i.i.d., and the number of which are approximately Poisson distributed. For a further discussion, see Aldous (1989).

In this paper, the main result is an upper bound for the total variation distance between the distribution of this quantity and a particular compound Poisson distribution. The bound can be expressed in terms of the first two moments of the interrenewal time conditional distributions, and on two constants obtained from each Radon–Nikodym derivative of an interrenewal time conditional distribution with respect to an exponential distribution, by solving a small number of systems of linear equations of dimension at most the total number of states. This is explicit often enough to be of considerable interest.

We briefly describe the ideas in the proof. If a single state $a \in B^c$ is chosen, we can construct a bound of the desired kind by expressing the quantity of interest as the accumulated reward of an embedded renewal reward process, for which the points with state $a$ serve as renewals. We then use Theorem 5.1 in Erhardsson (2000b) which gives a compound Poisson approximation error bound for the accumulated reward. However, the bound is small only if points with state $a$ are frequently observed. For many Markov chains, there exists a frequently observed state $a$ [see Erhardsson (1999, 2000a, 2001a, b)], but in many other cases no such $a$ exists.

To solve this problem, we study the pair of random variables $(\zeta, V)$, where $\zeta$ is the distance between two successive points and $V$ is the state of the second point. We construct a probability space containing $(\zeta, V)$ and a third random variable $\hat{\zeta}$ such that, for all $t$, conditional on $\{\hat{\zeta} \leq t, \zeta > t\}$, the pair $(\zeta - t, V)$ has the distribution $\nu_\gamma \times \mu$, where $\nu_\gamma$ is an exponential (or geometric) distribution with mean $\gamma^{-1}$, and $\mu$ is a fixed distribution. One might say that the event $\{\hat{\zeta} \leq t, \zeta > t\}$ indicates a loss of memory at or before $t$. For this reason, we call $\hat{\zeta}$ a "strong memoryless time."

Using strong memoryless times, we embed the stationary MRPP into another stationary MRPP whose state space contains an additional state 0. The points with states different from 0 also belong to the original MRPP. The points with state 0 represent losses of memory in the original MRPP, and are frequently observed if the original MRPP loses its memory quickly enough. The bound is then derived by an application of Theorem 5.1 in Erhardsson (2000b) to the accumulated reward of a renewal reward process embedded into the new MRPP, for which the points with state 0 serve as renewals.

In the last section, we compute the bound explicitly for an important special case: the number of points in $(0, t]$ of a stationary renewal process in continuous time. The bound is 0 if the interrenewal times are exponentially distributed, that is, if the renewal process is Poisson. We intend to present other applications of our results in the future.

It should be emphasized that the results in this paper are not limit theorems, but total variation distance error bounds which are valid for all finite parameter values. If desired, they can be used to derive limit theorems for various kinds of asymptotics, by showing that the bound converges to 0 under these asymptotics. They can also be used to bound the rate of convergence in limit theorems, by bounding the rate of convergence of the error bound.

It should also be mentioned that the literature contains a number of results concerning weak convergence to a compound Poisson point process, for special kinds of point processes (e.g., thinned point processes, or point processes generated by extreme values). Most of these are pure limit theorems without error bounds; see, for example, Serfozo (1984) and Leadbetter and Rootzén (1988). A few error bounds also exist, but not intended for processes of the kind studied in this paper, and derived using methods very different from ours; see, for example, Barbour and Månsson (2002).

The rest of the paper is organized as follows. In Section 2, some basic notation is given. In Section 3, we give necessary and sufficient conditions for the existence of strong memoryless times, and derive some of their relevant properties. In Section 4, we derive bounds for the total variation distance between the distribution of the number of points of an MRPP in $(0, t]$ with states in $B$ and a compound Poisson distribution. In Section 5, we consider the number of points in $(0, t]$ of a stationary renewal process, and obtain a more explicit expression for the bound.

**2. Basic notation.** Sets of numbers are denoted as follows: $\mathbb{R} =$ the real numbers, $\mathbb{Z} =$ the integers, $\mathbb{R}_+ = [0, \infty)$, $\mathbb{R}'_+ = (0, \infty)$, $\mathbb{Z}_+ = \{0, 1, 2, \ldots\}$ and $\mathbb{Z}'_+ = \{1, 2, \ldots\}$. The distribution of any random element $X$ in any measurable space $(S, \mathscr{S})$ is denoted by $\mathscr{L}(X)$. The Borel $\sigma$-algebra of any topological space $S$ is denoted by $\mathscr{B}_S$.

A compound Poisson distribution is a probability distribution with a characteristic function of the form $\phi(t) = \exp(-\int_{\mathbb{R}'_+} (1 - e^{itx}) \, d\pi(x))$, where $\pi$ is a measure on $(\mathbb{R}'_+, \mathscr{B}_{\mathbb{R}'_+})$ such that $\int_{\mathbb{R}'_+} (1 \wedge x) \, d\pi(x) < \infty$. It is denoted by $\text{POIS}(\pi)$. If $\|\pi\| = \pi(\mathbb{R}'_+) < \infty$, then $\text{POIS}(\pi) = \mathscr{L}(\sum_{i=1}^{U} T_i)$, where $\mathscr{L}(T_i) = \pi / \|\pi\|$ for each $i \in \mathbb{Z}'_+$, $U \sim \text{Po}(\|\pi\|)$, and all random variables are independent.

The total variation distance is a metric on the space of probability measures on any measurable space $(S, \mathscr{S})$. It is defined for two such measures

$\nu_1$ and $\nu_2$ by

$$d_{\mathrm{TV}}(\nu_1, \nu_2) = \sup_{A \in \mathscr{S}} |\nu_1(A) - \nu_2(A)|.$$

**3. Strong memoryless times.** In Theorems 3.1–3.3, we define strong memoryless times, give necessary and sufficient conditions for their existence, and derive some of their relevant properties. Note that Theorem 3.1 holds under more general conditions than are needed in Section 4. This will facilitate other applications in the future.

By $\nu_\gamma$ we mean the exponential distribution with mean $\gamma^{-1}$.

THEOREM 3.1. *Let $(\zeta, V)$ be a random variable taking values in $(\mathbb{R}_+ \times S, \mathscr{B}_{\mathbb{R}_+} \times \mathscr{S})$, where $(S, \mathscr{S})$ is a measurable space. Let $\mu$ be a probability measure on $(S, \mathscr{S})$. Assume that $\sigma: \mathbb{R}_+ \to [0, 1]$ satisfies*

$$(3.1) \qquad \sigma(t) \leq \inf_{\substack{C \in \mathscr{B}_{\mathbb{R}_+'} \times \mathscr{S} \\ (\nu_\gamma \times \mu)(C) > 0}} \frac{\mathbb{P}((\zeta - t, V) \in C)}{(\nu_\gamma \times \mu)(C)} \qquad \forall t \in \mathbb{R}_+,$$

*and that $G: \mathbb{R}_+ \to \mathbb{R}_+$, defined by $G(t) = \sigma(t)e^{\gamma t}$, is nondecreasing and right-continuous. In particular, these conditions are satisfied if equality holds in (3.1). Then we can define, on the same probability space as $(\zeta, V)$, a non-negative random variable $\hat{\zeta}$ (called a* strong memoryless time*) such that*

$$(3.2) \quad \begin{aligned} &\mathbb{P}(\hat{\zeta} \leq t, \zeta \leq u, V \in A) \\ &= \mathbb{P}(\zeta \leq t \wedge u, V \in A) + \sigma(t)(1 - e^{-\gamma[u-t]_+})\mu(A) \\ &\hspace{4cm} \forall (t, u, A) \in \mathbb{R}_+ \times \mathbb{R}_+ \times \mathscr{S}, \end{aligned}$$

*and such that $\mathbb{P}(\hat{\zeta} \leq \zeta) = 1$ and $\mathscr{L}((\zeta - t, V) | \hat{\zeta} \leq t, \zeta > t) = \nu_\gamma \times \mu$ for each $t \in \mathbb{R}_+$. Conversely, assume that the nonnegative random variable $\hat{\zeta}$, defined on the same probability space as $(\zeta, V)$, satisfies $\mathbb{P}(\hat{\zeta} \leq \zeta) = 1$ and $\mathscr{L}((\zeta - t, V) | \hat{\zeta} \leq t, \zeta > t) = \nu_\gamma \times \mu$ for each $t \in \mathbb{R}_+$. Then $\sigma: \mathbb{R}_+ \to [0, 1]$, defined by $\sigma(t) = \mathbb{P}(\hat{\zeta} \leq t, \zeta > t)$, satisfies (3.1), and $G: \mathbb{R}_+ \to \mathbb{R}_+$, defined by $G(t) = \sigma(t)e^{\gamma t}$, is nondecreasing and right-continuous.*

PROOF. For notational convenience, extend $\sigma$ to a function $\sigma: \mathbb{R} \to [0, 1]$ by defining $\sigma(t) = 0$ for each $t < 0$, and define $F: \mathbb{R} \times \mathbb{R} \times \mathscr{S} \to [0, 1]$ by $F(t, u, A) = \mathbb{P}(\zeta \leq t \wedge u, V \in A) + \sigma(t)(1 - e^{-\gamma[u-t]_+})\mu(A)$. It is easy to see that if we can define a random variable $(\hat{\zeta}, \zeta, V)$ taking values in $(\mathbb{R} \times \mathbb{R} \times S, \mathscr{B}_\mathbb{R} \times \mathscr{B}_\mathbb{R} \times \mathscr{S})$ such that

$$\mathbb{P}(\hat{\zeta} \leq t, \zeta \leq u, V \in A) = F(t, u, A) \qquad \forall (t, u, A) \in \mathbb{R} \times \mathbb{R} \times \mathscr{S},$$

then $\mathbb{P}(\hat{\zeta} \leq \zeta) = 1$ and $\mathscr{L}((\zeta - t, V)|\hat{\zeta} \leq t, \zeta > t) = \nu_\gamma \times \mu$ for each $t \in \mathbb{R}_+$. Hence, for the first part of the theorem it suffices to prove that there exists a probability distribution $\lambda_F$ on $(\mathbb{R} \times \mathbb{R} \times S, \mathscr{B}_\mathbb{R} \times \mathscr{B}_\mathbb{R} \times \mathscr{S})$ such that

$$\lambda_F((-\infty, t] \times (-\infty, u] \times A) = F(t, u, A) \qquad \forall (t, u, A) \in \mathbb{R} \times \mathbb{R} \times \mathscr{S}.$$

To do this, we use Theorem 11.3 in Billingsley (1986). Define $\mathscr{H}$ by

$$\mathscr{H} = \{(a,b] \times (c,d] \times A; -\infty < a \leq b < \infty, -\infty < c \leq d < \infty, A \in \mathscr{S}\}.$$

Clearly, $\mathscr{H}$ is a semiring generating $\mathscr{B}_\mathbb{R} \times \mathscr{B}_\mathbb{R} \times \mathscr{S}$. Define a set function $\lambda_F : \mathscr{H} \to \mathbb{R}$ by

$$\begin{aligned}
&\lambda_F((a,b] \times (c,d] \times A) \\
&= F(b,d,A) - F(a,d,A) - F(b,c,A) + F(a,c,A) \\
&= \mathbb{P}(\zeta \in (a,b] \cap (c,d], V \in A) + \sigma(b)(e^{-\gamma[c-b]_+} - e^{-\gamma[d-b]_+})\mu(A) \\
&\quad - \sigma(a)(e^{-\gamma[c-a]_+} - e^{-\gamma[d-a]_+})\mu(A) \qquad \forall (a,b] \times (c,d] \times A \in \mathscr{H}.
\end{aligned}$$

Using the facts that $\sigma$ satisfies (3.1) and that $G$ is nondecreasing, it can be shown that $\lambda_F$ is nonnegative. For example, if $a < b \leq c < d$, we get

$$\begin{aligned}
&\lambda_F((a,b] \times (c,d] \times A) \\
&= \sigma(b)(e^{-\gamma(c-b)} - e^{-\gamma(d-b)})\mu(A) \\
&\quad - \sigma(a)(e^{-\gamma(c-a)} - e^{-\gamma(d-a)})\mu(A) \\
&= (e^{-\gamma c} - e^{-\gamma d})(\sigma(b)e^{\gamma b} - \sigma(a)e^{\gamma a})\mu(A) \geq 0,
\end{aligned}$$

while if $a \leq c < b \leq d$ we get

$$\begin{aligned}
&\lambda_F((a,b] \times (c,d] \times A) \\
&= \mathbb{P}(\zeta \in (c,b], V \in A) + \sigma(b)(1 - e^{-\gamma(d-b)})\mu(A) \\
&\quad - \sigma(a)(e^{-\gamma(c-a)} - e^{-\gamma(d-a)})\mu(A) \\
&= \mathbb{P}(\zeta \in (c,b], V \in A) + \sigma(b)e^{\gamma b}(e^{-\gamma b} - e^{-\gamma d})\mu(A) \\
&\quad - \sigma(a)e^{\gamma a}(e^{-\gamma c} - e^{-\gamma d})\mu(A) \\
&\geq \mathbb{P}(\zeta \in (c,b], V \in A) - \sigma(a)e^{\gamma a}(e^{-\gamma c} - e^{-\gamma b})\mu(A) \geq 0.
\end{aligned}$$

We now show that $\lambda_F$ is countably additive on $\mathscr{H}$. In other words, we assume that $(a,b] \times (c,d] \times A = \bigcup_{i=1}^\infty (a_i, b_i] \times (c_i, d_i] \times A_i$, where $(a,b] \times (c,d] \times A \in \mathscr{H}$, $(a_i, b_i] \times (c_i, d_i] \times A_i \in \mathscr{H}$ for each $i \in \mathbb{Z}'_+$, and the sets $\{(a_i, b_i] \times (c_i, d_i] \times A_i; i \in \mathbb{Z}'_+\}$ are disjoint, and show that

(3.3) $$\lambda_F((a,b] \times (c,d] \times A) = \sum_{i=1}^\infty \lambda_F((a_i, b_i] \times (c_i, d_i] \times A_i).$$

Define $F_A : \mathbb{R} \times \mathbb{R} \to [0,1]$ by $F_A(t,u) = F(t,u,A)$ [where $A$ is the same set as in (3.3)]. Define also the semiring $\mathscr{H}^*$ and the set function $\lambda_{F_A} : \mathscr{H}^* \to \mathbb{R}$ by

$$\mathscr{H}^* = \{(a,b] \times (c,d]; -\infty < a \leq b < \infty, -\infty < c \leq d < \infty\};$$
$$\lambda_{F_A}((a,b] \times (c,d]) = \lambda_F((a,b] \times (c,d] \times A) \qquad \forall\, (a,b] \times (c,d] \in \mathscr{H}^*.$$

Clearly, $F_A$ is continuous from above, and it was shown earlier that $\lambda_{F_A}$ is nonnegative. It therefore follows from Theorem 12.5 in Billingsley (1986) that $\lambda_{F_A}$ can be uniquely extended to a measure on $(\mathbb{R} \times \mathbb{R}, \mathscr{B}_\mathbb{R} \times \mathscr{B}_\mathbb{R})$, which in turn implies that $\lambda_{F_A} \times \mu$ is a measure on $(\mathbb{R} \times \mathbb{R} \times S, \mathscr{B}_\mathbb{R} \times \mathscr{B}_\mathbb{R} \times \mathscr{S})$. Hence,

$$\lambda_{F_A}((a,b] \times (c,d])\mu(A) = \sum_{i=1}^{\infty} \lambda_{F_A}((a_i,b_i] \times (c_i,d_i])\mu(A_i),$$

from which (3.3) will follow if we can show that

$$\sum_{i=1}^{\infty} \mathbb{P}(\zeta \in (a_i,b_i] \cap (c_i,d_i], V \in A)\mu(A_i)$$
$$= \sum_{i=1}^{\infty} \mathbb{P}(\zeta \in (a_i,b_i] \cap (c_i,d_i], V \in A_i)\mu(A).$$

But this follows from the facts that

$$\mathbb{P}(\zeta \in (a,b] \cap (c,d], V \in A)\mu(A) = \sum_{i=1}^{\infty} \mathbb{P}(\zeta \in (a_i,b_i] \cap (c_i,d_i], V \in A)\mu(A_i)$$

and

$$\mathbb{P}(\zeta \in (a,b] \cap (c,d], V \in A) = \sum_{i=1}^{\infty} \mathbb{P}(\zeta \in (a_i,b_i] \cap (c_i,d_i], V \in A_i).$$

This concludes the proof of the first part of the theorem.

We next show that if $\sigma$ is chosen so that equality holds in (3.1), then $G$ is nondecreasing and right-continuous. Let $C \in \mathscr{B}_{\mathbb{R}'_+} \times \mathscr{S}$ and define, for each $t \in \mathbb{R}_+$, $C^t = \{(x+t, y); (x,y) \in C\}$. It is easy to show that $(\nu_\gamma \times \mu)(C^t) = e^{-\gamma t}(\nu_\gamma \times \mu)(C)$ for each $t \in \mathbb{R}_+$. Hence, for each $0 \leq s < t < \infty$,

$$\frac{\mathbb{P}((\zeta - t, V) \in C)}{(\nu_\gamma \times \mu)(C)} = \frac{\mathbb{P}((\zeta - s, V) \in C^{t-s})e^{-\gamma(t-s)}}{(\nu_\gamma \times \mu)(C^{t-s})},$$

implying that

$$G(t) = \inf_{\substack{C \in \mathscr{B}_{\mathbb{R}'_+} \times \mathscr{S} \\ (\nu_\gamma \times \mu)(C) > 0}} \frac{\mathbb{P}((\zeta - t, V) \in C)e^{\gamma t}}{(\nu_\gamma \times \mu)(C)}$$

$$= \inf_{\substack{C \in \mathscr{B}_{\mathbb{R}'_+} \times \mathscr{S} \\ (\nu_\gamma \times \mu)(C) > 0}} \frac{\mathbb{P}((\zeta - s, V) \in C^{t-s})e^{\gamma s}}{(\nu_\gamma \times \mu)(C^{t-s})} \geq G(s),$$

so $G$ is nondecreasing. Next, fix $t \in \mathbb{R}_+$ and choose a sequence $\{C_k \in \mathscr{B}_{\mathbb{R}'_+} \times \mathscr{S}; k \in \mathbb{Z}'_+\}$ such that $(\nu_\gamma \times \mu)(C_k) > 0$ for each $k \in \mathbb{Z}'_+$ and

$$\lim_{k \to \infty} \frac{\mathbb{P}((\zeta - t, V) \in C_k)}{(\nu_\gamma \times \mu)(C_k)} = \inf_{\substack{C \in \mathscr{B}_{\mathbb{R}'_+} \times \mathscr{S} \\ (\nu_\gamma \times \mu)(C) > 0}} \frac{\mathbb{P}((\zeta - t, V) \in C)}{(\nu_\gamma \times \mu)(C)}.$$

For each $k \in \mathbb{Z}'_+$ and $u \in \mathbb{R}_+$, define $C_{k,u} = C_k \cap ((u, \infty) \times S)$ and $C_{k,u}^{-u} = \{(x - u, y); (x, y) \in C_{k,u}\}$. Then, for each $k \in \mathbb{Z}'_+$ and each $u \in \mathbb{R}_+$ such that $(\nu_\gamma \times \mu)(C_{k,u}) > 0$,

$$G(t+u) = \inf_{\substack{C \in \mathscr{B}_{\mathbb{R}'_+} \times \mathscr{S} \\ (\nu_\gamma \times \mu)(C) > 0}} \frac{\mathbb{P}((\zeta - t - u, V) \in C)e^{\gamma(t+u)}}{(\nu_\gamma \times \mu)(C)}$$

$$\leq \frac{\mathbb{P}((\zeta - t - u, V) \in C_{k,u}^{-u})e^{\gamma(t+u)}}{(\nu_\gamma \times \mu)(C_{k,u}^{-u})} = \frac{\mathbb{P}((\zeta - t, V) \in C_{k,u})e^{\gamma t}}{(\nu_\gamma \times \mu)(C_{k,u})}.$$

This implies that $\limsup_{u \downarrow 0} G(t + u) \leq G(t)$, and since $G$ is nondecreasing, it must be right-continuous.

For the last part of the theorem, assume that a nonnegative random variable $\hat{\zeta}$ can be defined on the same probability space as $(\zeta, V)$, such that $\mathbb{P}(\hat{\zeta} \leq \zeta) = 1$ and $\mathscr{L}((\zeta - t, V)|\hat{\zeta} \leq t, \zeta > t) = \nu_\gamma \times \mu$ for each $t \in \mathbb{R}_+$. Then,

$$\mathbb{P}((\zeta - t, V) \in C) = (\nu_\gamma \times \mu)(C)\mathbb{P}(\hat{\zeta} \leq t, \zeta > t)$$
$$+ \mathbb{P}((\zeta - t, V) \in C|\hat{\zeta} > t, \zeta > t)\mathbb{P}(\hat{\zeta} > t, \zeta > t)$$
$$\forall (t, C) \in \mathbb{R}_+ \times (\mathscr{B}_{\mathbb{R}'_+} \times \mathscr{S}),$$

which implies that $\sigma : \mathbb{R}_+ \to [0, 1]$, defined by $\sigma(t) = \mathbb{P}(\hat{\zeta} \leq t, \zeta > t)$, satisfies (3.1). Moreover, (3.2) holds with $\sigma$ defined in this way, which implies that if $a < b \leq c < d$, then

$$\mathbb{P}(\hat{\zeta} \in (a, b], \zeta \in (c, d])$$
$$= \sigma(b)(e^{-\gamma(c-b)} - e^{-\gamma(d-b)}) - \sigma(a)(e^{-\gamma(c-a)} - e^{-\gamma(d-a)})$$
$$= (e^{-\gamma c} - e^{-\gamma d})(\sigma(b)e^{\gamma b} - \sigma(a)e^{\gamma a}) \geq 0,$$

so $G$ is nondecreasing, and clearly also right-continuous. □

THEOREM 3.2. *Let $(\hat{\zeta}, \zeta, V)$ be a random variable taking values in $(\mathbb{R}_+ \times \mathbb{R}_+ \times S, \mathscr{B}_{\mathbb{R}_+} \times \mathscr{B}_{\mathbb{R}_+} \times \mathscr{S})$, where $(S, \mathscr{S})$ is a measurable space. Let $\mu$ be a* probability measure on $(S, \mathscr{S})$. Define $\sigma: \mathbb{R}_+ \to [0,1]$ by $\sigma(t) = \mathbb{P}(\hat{\zeta} \leq t, \zeta > t)$. If $\mathbb{P}(\hat{\zeta} \leq \zeta) = 1$ and $\mathscr{L}((\zeta - t, V) | \hat{\zeta} \leq t, \zeta > t) = \nu_\gamma \times \mu$ for each $t \in \mathbb{R}_+$, then $\mathscr{L}((\hat{\zeta}, \zeta - \hat{\zeta}, V) | \hat{\zeta} < \zeta) = \mathscr{L}(\hat{\zeta} | \hat{\zeta} < \zeta) \times \nu_\gamma \times \mu$, where

$$\mathbb{P}(\hat{\zeta} \leq t, \hat{\zeta} < \zeta) = \sigma(t) + \gamma \int_0^t \sigma(x)\, dx \qquad \forall t \in \mathbb{R}_+$$

and

$$\mathbb{P}(\hat{\zeta} = \zeta \leq t, V \in A)$$
$$= \mathbb{P}(\zeta \leq t, V \in A) - \mu(A)\gamma \int_0^t \sigma(x)\, dx \qquad \forall (t, A) \in \mathbb{R}_+ \times \mathscr{S}.$$

Conversely, if $\mathbb{P}(\hat{\zeta} \leq \zeta) = 1$ and $\mathscr{L}((\hat{\zeta}, \zeta - \hat{\zeta}, V) | \hat{\zeta} < \zeta) = \mathscr{L}(\hat{\zeta} | \hat{\zeta} < \zeta) \times \nu_\gamma \times \mu$, then $\mathscr{L}((\zeta - t, V) | \hat{\zeta} \leq t, \zeta > t) = \nu_\gamma \times \mu$ for each $t \in \mathbb{R}_+$.

PROOF. From (3.2), and using bounded convergence, we get

$$\mathbb{P}(\hat{\zeta} \in (0, t], \zeta - \hat{\zeta} \in (0, u], V \in A)$$
$$= \lim_{N \to \infty} \sum_{i=1}^N \mathbb{P}\left(\hat{\zeta} \in \left(\frac{(i-1)t}{N}, \frac{it}{N}\right], \zeta \in \left(\frac{it}{N}, \frac{it}{N} + u\right], V \in A\right)$$
$$= \mu(A) \lim_{N \to \infty} \left((1 - e^{-\gamma u}) \sum_{i=1}^N \sigma\left(\frac{it}{N}\right)\right.$$
$$\left. - (e^{-\gamma(t/N)} - e^{-\gamma(t/N + u)}) \sum_{i=1}^N \sigma\left(\frac{(i-1)t}{N}\right)\right)$$
$$= \mu(A)(1 - e^{-\gamma u}) \lim_{N \to \infty} \sum_{i=1}^N \left(\sigma\left(\frac{it}{N}\right) - \sigma\left(\frac{(i-1)t}{N}\right)\right)$$
$$+ \mu(A)(1 - e^{-\gamma u}) \lim_{N \to \infty} (1 - e^{-\gamma(t/N)}) \sum_{i=1}^N \sigma\left(\frac{(i-1)t}{N}\right)$$
$$\forall (t, u, A) \in \mathbb{R}_+ \times \mathbb{R}_+ \times \mathscr{S}.$$

The first sum telescopes. For the second sum, we note that $\sigma$ is Riemann integrable on $[0, t]$. This holds since the function $G: \mathbb{R}_+ \to \mathbb{R}_+$, defined by $G(t) = \sigma(t)e^{\gamma t}$, is nondecreasing, hence Riemann–Stieltjes integrable on $[0, t]$ with respect to $\alpha: \mathbb{R}_+ \to [0, 1]$, defined by $\alpha(t) = 1 - e^{-\gamma t}$; see Theorems 6.9 and 6.17 in Rudin (1976). This gives

$$\mathbb{P}(\hat{\zeta} \in (0, t], \zeta - \hat{\zeta} \in (0, u], V \in A)$$
$$= \mu(A)(1 - e^{-\gamma u})\left(\sigma(t) - \sigma(0) + \gamma \int_0^t \sigma(x)\, dx\right)$$

$$\forall (t, u, A) \in \mathbb{R}_+ \times \mathbb{R}_+ \times \mathscr{S}.$$

To complete the proof of the first part of the theorem, note that

$$\mathbb{P}(\hat{\zeta} = 0, \zeta \in (0, u], V \in A) = \mu(A)(1 - e^{-\gamma u})\sigma(0),$$

and that $\mathbb{P}(\hat{\zeta} = \zeta \leq t, V \in A) = \mathbb{P}(\hat{\zeta} \leq t, V \in A) - \mathbb{P}(\hat{\zeta} \leq t, \hat{\zeta} < \zeta, V \in A)$. For the second part of the theorem,

$$\begin{aligned}
\mathbb{P}(\zeta - t &\leq u, V \in A, \hat{\zeta} \leq t, \zeta > t) \\
&= \mathbb{E}(e^{\gamma \hat{\zeta}} I\{\hat{\zeta} \leq t, \hat{\zeta} < \zeta\})(e^{-\gamma t} - e^{-\gamma(t+u)})\mu(A) \\
&= \mathbb{P}(\hat{\zeta} \leq t, \zeta > t)(1 - e^{-\gamma u})\mu(A) \qquad \forall (t, u, A) \in \mathbb{R}_+ \times \mathbb{R}_+ \times \mathscr{S}. \quad \square
\end{aligned}$$

THEOREM 3.3. *Let the conditions of Theorem 3.1 hold with $S$ a finite set, and let $f : \mathbb{R}_+ \times S \to \mathbb{R}_+$ be the Radon–Nikodym derivative with respect to $\nu_\gamma \times \mu$ of the part of $\mathscr{L}(\zeta, V)$ which is absolutely continuous with respect to $\nu_\gamma \times \mu$. Then,*

$$\inf_{\substack{C \in \mathscr{B}_{\mathbb{R}'_+} \times \mathscr{S} \\ (\nu_\gamma \times \mu)(C) > 0}} \frac{\mathbb{P}((\zeta - t, V) \in C)}{(\nu_\gamma \times \mu)(C)} = e^{-\gamma t} \operatorname*{ess\,inf}_{x \in (t, \infty) \times S} f(x) \qquad \forall t \in \mathbb{R}_+.$$

PROOF. The "$\geq$" part is easy. For the "$\leq$" part, we use Theorem 35.8 in Billingsley (1986). For each $n \in \mathbb{Z}_+$, let $\mathscr{F}_n$ be the $\sigma$-algebra generated by the sets $\{(k2^{-n}, (k+1)2^{-n}] \times \{s\}; k \in \mathbb{Z}_+, s \in S\}$. It is well known that $\sigma(\bigcup_{n=0}^\infty \mathscr{F}_n) = \mathscr{B}_{\mathbb{R}'_+} \times \mathscr{B}_S$. Therefore, for $\nu_\gamma \times \mu$–almost every $x \in \mathbb{R}_+ \times S$, $f(x)$ is the limit of ratios of the kind appearing on the left-hand side. $\quad\square$

REMARK 3.1. The strong memoryless time $\hat{\zeta}$ for which equality holds in (3.1) is optimal in the sense that $\mathbb{P}(\hat{\zeta} \leq t | \zeta > t)$ is maximized uniformly over all $t \in \mathbb{R}_+$.

REMARK 3.2. Theorems 3.1 and 3.2 imply that $\hat{\zeta}$ is a strong memoryless time for $(\zeta, V)$ if and only if $(\hat{\zeta}, \zeta - \hat{\zeta}, V) = \chi(\eta_0, \eta_1, V_1) + (1 - \chi)(\eta_2, 0, V_2)$, where the random variables $\chi, \eta_0, \eta_1, V_1$ and $(\eta_2, V_2)$ are independent, $\chi$ takes values in $\{0, 1\}$, $\eta_1$ is exponentially distributed with mean $\gamma^{-1}$ and $\mathscr{L}(V_1) = \mu$. Clearly, $\sigma(t) = \mathbb{P}(\hat{\zeta} \leq t, \zeta > t) = \mathbb{P}(\chi = 1)\mathbb{E}(e^{-\gamma(t - \eta_0)} I\{\eta_0 \leq t\})$.

REMARK 3.3. Let $S = \{1\}$, and let $f$ be the Radon–Nikodym derivative of $\mathscr{L}(\zeta)$ with respect to the exponential distribution with mean $\gamma^{-1}$.

1. Assume that $f(t) \geq \lim_{u \to \infty} f(u) = c > 0$ for all $t \in \mathbb{R}_+$. Then, the optimal choice of $\sigma$ is $\sigma(t) = ce^{-\gamma t}$ which, by Theorem 3.2, implies that $\mathbb{P}(\chi = 1) = c$ and $\eta_0 \equiv 0$.

2. Assume that $f$ is nondecreasing. Then, the optimal choice of $\sigma$ is $\sigma(t) = f(t)e^{-\gamma t}$ which, again by Theorem 3.2, implies that $\chi \equiv 1$ and $\mathbb{P}(\eta_0 \leq t) = f(t)e^{-\gamma t} + \mathbb{P}(\zeta \leq t)$ for each $t \in \mathbb{R}_+$.

REMARK 3.4. The strong memoryless times were originally inspired by another construction, the strong stationary times used in Aldous and Diaconis (1986, 1987) and Diaconis and Fill (1990) to bound the rate of convergence of a finite-state discrete-time Markov chain $\{\eta_i; i \in \mathbb{Z}_+\}$ to the stationary distribution $\mu$. A strong stationary time $T$ is a randomized stopping time such that $\mathscr{L}(\eta_i|T \leq i) = \mu$ for each $i \in \mathbb{Z}_+$. It seems unlikely that strong stationary times could be used (even in the restricted setting of discrete-time Markov chains) to solve the problem considered in the present paper, without significant modifications leading in the end to the construction of strong memoryless times.

Strong memoryless times are also related to a construction due to Athreya and Ney (1978) and Nummelin (1978), known as splitting. This is an embedding of a discrete-time Markov chain on a general state space (satisfying an irreducibility condition) into another Markov chain on a larger state space which contains a recurrent single state. In general, this recurrent state need not be frequently observed, so splitting does not suffice (even in the discrete-time Markov chain setting) to solve the problem considered in the present paper.

We end this section with lattice versions of the preceding theorems. The proofs are analogous to those above, but simpler, since right-continuity is trivial in the lattice case.

THEOREM 3.4. *Let the conditions of Theorem 3.1 hold, with the following changes*: $\mathbb{R}_+$ *is replaced by* $\mathbb{Z}_+$, $\nu_\gamma$ *is the geometric distribution with mean* $\gamma^{-1}$, *and* $e^{-\gamma}$ *is replaced by* $1 - \gamma$ *in the definition of $G$ and in* (3.2). *Then, all the assertions of Theorem 3.1 remain valid.*

THEOREM 3.5. *Let the conditions of Theorem 3.2 hold, with the following changes*: $\mathbb{R}_+$ *is replaced by* $\mathbb{Z}_+$, *and* $\nu_\gamma$ *is the geometric distribution with mean* $\gamma^{-1}$. *Then, all the assertions of Theorem 3.2 remain valid, with* $\int_0^t \sigma(x) \, dx$ *replaced by* $\sum_{i=0}^{t-1} \sigma(i)$.

**4. Markov renewal point processes.** In this section we use the results in Section 3 to address the problem described in Section 1. Recall that we wish to find a bound for the total variation distance between the distribution of the number of points of an MRPP in $(0, t]$ with states in $B$, and a suitable compound Poisson distribution. We assume that the reader is familiar with the basic theory of marked point processes. Good references are Rolski (1981), Franken, König, Arndt and Schmidt (1982) and Port and Stone (1973).

We begin with the definition of an MRPP. Let $S = \{1, \ldots, N\}$, and let $\{(\zeta_i^S, V_{i+1}^S); i \in \mathbb{Z}\}$ be a stationary discrete-time Markov chain taking values in $(\mathbb{R}_+ \times S, \mathscr{B}_{\mathbb{R}_+ \times S})$, with a transition probability $p$ such that $p((t,s), \cdot) = p(s, \cdot)$ for each $(t,s) \in \mathbb{R}_+ \times S$. Assume that $\{V_i^S; i \in \mathbb{Z}\}$ is irreducible, and that $0 < \mathbb{E}(\zeta_0^S) < \infty$. (We collectively denote these conditions by C0.)

For each $A \subset S$, let $\{(\zeta_i^A, V_{i+1}^A); i \in \mathbb{Z}\}$ have the distribution $\mathscr{L}((\zeta_i^S, V_{i+1}^S); i \in \mathbb{Z} | V_0^S \in A)$, and define $\{U_i^A; i \in \mathbb{Z}\}$ by $U_0^A = 0$, $U_i^A = \sum_{j=0}^{i-1} \zeta_j^A$ for each $i \geq 1$, and $U_i^A = -\sum_{j=i}^{-1} \zeta_j^A$ for each $i \leq -1$. Define the point process $\Psi^A$ on $(\mathbb{R} \times S, \mathscr{B}_{\mathbb{R} \times S})$ by $\Psi^A(\cdot) = \sum_{i \in \mathbb{Z}} I\{(U_i^A, V_i^A) \in \cdot\}$. $\Psi^A$ is a Palm version (with respect to marks in $A$) of an MRPP.

Next, define the point process $\Psi$ on $(\mathbb{R} \times S, \mathscr{B}_{\mathbb{R} \times S})$ by

$$(4.1) \qquad \mathbb{E}(g(\Psi)) = \frac{\mathbb{E}(\int_0^{U_{\tau_1^A}^A} g(\theta_t(\Psi^A)) \, dt)}{\mathbb{E}(U_{\tau_1^A}^A)} \qquad \forall g \in \mathcal{F}_{\mathcal{N}(\mathbb{R} \times S)}^+,$$

where $\mathcal{F}_{\mathcal{N}(\mathbb{R} \times S)}^+$ are the nonnegative Borel functions on the space of counting measures on $(\mathbb{R} \times S, \mathscr{B}_{\mathbb{R} \times S})$, $\tau_1^A = \min\{i \geq 1; V_i^A \in A\}$ and $\theta$ is the shift operator, defined by $\theta_t(\Psi)((a,b] \times \cdot) = \Psi((a+t, b+t] \times \cdot)$. This definition is independent of the choice of $A$, and $\Psi$ is a stationary marked point process. There exist random variables $\{(U_i, V_i); i \in \mathbb{Z}\}$ (where $\cdots \leq U_{-1} \leq U_0 \leq 0 < U_1 \leq \cdots$) such that $\Psi(\cdot) = \sum_{i \in \mathbb{Z}} I\{(U_i, V_i) \in \cdot\}$. $\Psi$ is a stationary MRPP.

The quantity that we are interested in can be expressed as $\Psi((0, t] \times B)$. We assume without loss of generality that $B = S$, since otherwise we can replace $\Psi$ by its restriction to $(\mathbb{R} \times B, \mathscr{B}_{\mathbb{R} \times B})$, which is also a stationary MRPP.

Analogously, we may define, using a stationary discrete-time Markov chain $\{(\zeta_i^S, V_{i+1}^S); i \in \mathbb{Z}\}$ taking values in $(\mathbb{Z}_+ \times S, \mathscr{B}_{\mathbb{Z}_+ \times S})$, a stationary MRPP in discrete time. In this case, for each $A \subset S$, the distribution of $\Psi$ is given by a discrete version of (4.1), where the integral is replaced by a sum over the integers $\{0, \ldots, U_{\tau_1^A}^A - 1\}$.

We now explain how to use strong memoryless times to embed a stationary MRPP into another stationary MRPP which has favorable properties from the point of view of compound Poisson approximation. Consider a stationary discrete-time Markov chain $\{(\zeta_i^S, V_{i+1}^S); i \in \mathbb{Z}\}$ on the state space $(\mathbb{R}_+ \times S, \mathscr{B}_{\mathbb{R}_+ \times S})$ with transition probability $p$, satisfying condition C0. Denote by $\nu_\gamma$ the exponential distribution with mean $\gamma^{-1}$, and let $\mu$ be a probability measure on $(S, \mathscr{B}_S)$. For each $s \in S$, assume that $\sigma_s : \mathbb{R}_+ \to [0,1]$

satisfies

$$(4.2) \quad \sigma_s(t) \leq \inf_{\substack{C \in \mathscr{B}_{\mathbb{R}'_+} \times \mathscr{B}_S \\ (\nu_\gamma \times \mu)(C) > 0}} \frac{\mathbb{P}((\zeta_0^S - t, V_1^S) \in C | V_0^S = s)}{(\nu_\gamma \times \mu)(C)} \quad \forall\, t \in \mathbb{R}_+,$$

and that $G_s : \mathbb{R}_+ \to \mathbb{R}_+$, defined by $G_s(t) = \sigma_s(t) e^{\gamma t}$, is nondecreasing and right-continuous. Assume also that $\int_0^\infty \sigma_s(t)\, dt > 0$ for at least one $s \in S$. (We collectively denote these conditions by C1.) Let $\widetilde{S} = S \cup \{0\}$, and let $\{(\widetilde{\zeta}_i^{\widetilde{S}}, \widetilde{V}_{i+1}^{\widetilde{S}}); i \in \mathbb{Z}\}$ be a stationary discrete-time Markov chain on the state space $(\mathbb{R}_+ \times \widetilde{S}, \mathscr{B}_{\mathbb{R}_+ \times \widetilde{S}})$, with a transition probability $\widetilde{p}$ defined for each $(s, s') \in S \times S$ by

$$\widetilde{p}(s, [0, u] \times \{0\}) = \sigma_s(u) + \gamma \int_0^u \sigma_s(t)\, dt,$$

$$\widetilde{p}(s, [0, u] \times \{s'\}) = p(s, [0, u] \times \{s'\}) - \mu(s') \gamma \int_0^u \sigma_s(t)\, dt,$$

$$\widetilde{p}(0, [0, u] \times \{0\}) = (1 - \varepsilon)(1 - e^{-(\gamma/\varepsilon)u}),$$

$$\widetilde{p}(0, [0, u] \times \{s'\}) = \mu(s') \varepsilon (1 - e^{-(\gamma/\varepsilon)u}),$$

where $\varepsilon \in (0, 1)$. For each $A \subset \widetilde{S}$, let $\widetilde{\Psi}^A(\cdot) = \sum_{i \in \mathbb{Z}} I\{(\widetilde{U}_i^A, \widetilde{V}_i^A) \in \cdot\}$ be the Palm version (with respect to marks in $A$) of the MRPP associated with $\{(\widetilde{\zeta}_i^{\widetilde{S}}, \widetilde{V}_{i+1}^{\widetilde{S}}); i \in \mathbb{Z}\}$, and let $\Psi^0 = \Psi^{\{0\}}$. $\widetilde{\Psi}^A$ is a point process on $(\mathbb{R} \times \widetilde{S}, \mathscr{B}_{\mathbb{R} \times \widetilde{S}})$. Heuristically, 0 is a frequently observed state for $\widetilde{\Psi}^A$ if $\varepsilon$ is small enough, and if the MRPP $\Psi^A$ loses its memory quickly enough after each occurrence of a point. Let also $\widetilde{\Psi}(\cdot) = \sum_{i \in \mathbb{Z}} I\{(\widetilde{U}_i, \widetilde{V}_i) \in \cdot\}$ be the stationary MRPP associated with $\{(\widetilde{\zeta}_i^{\widetilde{S}}, \widetilde{V}_{i+1}^{\widetilde{S}}); i \in \mathbb{Z}\}$.

The following fact is now crucial, since it implies that we have constructed an embedding: the restriction of $\widetilde{\Psi}$ to $(\mathbb{R} \times S, \mathscr{B}_{\mathbb{R} \times S})$ has the same distribution as $\Psi$. To see this, let $C_1, \ldots, C_k$ be disjoint subsets of $\mathbb{R} \times S$, let $C_i^t = \{(x + t, y); (x, y) \in C_i\}$ for each $t \in \mathbb{R}_+$ and let $n_1, \ldots, n_k$ be nonnegative integers. Applying (4.1) with $A = S$ gives

$$\mathbb{E}\left(\prod_{i=1}^k I\{\widetilde{\Psi}(C_i) = n_i\}\right) = \frac{\mathbb{E}(\int_0^{\widetilde{U}_{\tau_1^S}^S} \prod_{i=1}^k I\{\widetilde{\Psi}^S(C_i^t) = n_i\}\, dt)}{\mathbb{E}(\widetilde{U}_{\tau_1^S}^S)}.$$

Clearly, we may replace $\widetilde{\Psi}^S(\cdot)$ by $\sum_{i \in \mathbb{Z}} I\{(\widetilde{U}_{\tau_i^S}^S, \widetilde{V}_{\tau_i^S}^S) \in \cdot\}$, where $\cdots \leq \tau_{-1}^S \leq \tau_0^S = 0 < \tau_1^S \leq \cdots$ are the random integers $\{i \in \mathbb{Z}; \widetilde{V}_i^S \in S\}$. It is straightforward to show, using Theorems 3.1 and 3.2 and the strong Markov property, that the random sequence $\{(\widetilde{U}_{\tau_{i+1}^S}^S - \widetilde{U}_{\tau_i^S}^S, \widetilde{V}_{\tau_{i+1}^S}^S); i \in \mathbb{Z}\}$ is a stationary Markov chain with transition probability $p$, that is, it has the same distribution as $\{(\zeta_i^S, V_{i+1}^S); i \in \mathbb{Z}\}$. Hence, $\{(\widetilde{U}_{\tau_i^S}^S, \widetilde{V}_{\tau_i^S}^S); i \in \mathbb{Z}\}$ has the same distribution as $\{(U_i^S, V_i^S); i \in \mathbb{Z}\}$, and since $\Psi^S(\cdot) = \sum_{i \in \mathbb{Z}} I\{(U_i^S, V_i^S) \in \cdot\}$, the proof is complete.

Finally, we need the following tools. Define $\{(X_i^0, Y_i^0); i \in \mathbb{Z}\}$ by $(X_i^0, Y_i^0) = (\widetilde{U}_{\tau_i^0}^0, \tau_{i+1}^0 - \tau_i^0 - 1)$, where $\cdots < \tau_{-1}^0 < \tau_0^0 = 0 < \tau_1^0 < \cdots$ are the random integers $\{i \in \mathbb{Z}; \widetilde{V}_i^0 = 0\}$. The strong Markov property implies that $\{(X_{i+1}^0 - X_i^0, Y_i^0); i \in \mathbb{Z}\}$ is an i.i.d. sequence. Let $\xi^0(\cdot) = \sum_{i \in \mathbb{Z}} I\{(X_i^0, Y_i^0) \in \cdot\}$ be a point process on $(\mathbb{R} \times \mathbb{Z}_+, \mathscr{B}_{\mathbb{R} \times \mathbb{Z}_+})$. By definition, this is a Palm version of a renewal reward process. Similarly, define $\{(X_i, Y_i); i \in \mathbb{Z}\}$ by $(X_i, Y_i) = (\widetilde{U}_{\tau_i}, \tau_{i+1} - \tau_i - 1)$, where $\cdots < \tau_{-1} < \tau_0 \leq 0 < \tau_1 < \cdots$ are the random integers $\{i \in \mathbb{Z}; \widetilde{V}_i = 0\}$, and let $\xi(\cdot) = \sum_{i \in \mathbb{Z}} I\{(X_i, Y_i) \in \cdot\}$ be a point process on $(\mathbb{R} \times \mathbb{Z}_+, \mathscr{B}_{\mathbb{R} \times \mathbb{Z}_+})$. It is straightforward to show that $\xi$ is the stationary renewal reward process corresponding to $\xi^0$.

It is now easy to state and prove the main result of this section. It will be demonstrated in Section 5 that the bound given below can be expressed in terms of a small number of parameters obtained from the functions $\{\sigma_s; s \in S\}$, by solving a small number of systems of linear equations.

THEOREM 4.1. *Let $\Psi$ be a stationary MRPP with state space $S = \{1, \ldots, N\}$, satisfying condition* C0 *above. Let $\gamma > 0$, let $\mu$ be a probability measure on $(S, \mathscr{B}_S)$ and assume that the functions $\{\sigma_s : \mathbb{R}_+ \to [0,1]; s \in S\}$ satisfy condition* C1 *above. Then,*

$$\begin{aligned}(4.3)\quad & d_{\mathrm{TV}}(\mathscr{L}(\Psi((0,t] \times S)), \mathrm{POIS}(\pi)) \\ & \leq \frac{2\mathbb{E}(\widetilde{U}_{\tau_1^0-1}^0)}{\mathbb{E}(X_1^0)} + H_1(\pi) \frac{3t\mathbb{E}(Y_0^0)}{\mathbb{E}(X_1^0)} \left( \frac{\mathbb{E}(X_1^0 Y_0^0)}{\mathbb{E}(X_1^0)} + \frac{\mathbb{E}((X_1^0)^2)\mathbb{E}(Y_0^0)}{\mathbb{E}(X_1^0)^2} \right), \end{aligned}$$

*where $\pi_i = \frac{t}{\mathbb{E}(X_1^0)} \mathbb{P}(Y_0^0 = i)$ for $i \geq 1$, and*

$$H_1(\pi) \leq \begin{cases} \left(1 \wedge \frac{1}{\pi_1}\right) e^{\|\pi\|}, & \textit{always,} \\ 1 \wedge \frac{1}{\pi_1 - 2\pi_2} \left( \frac{1}{4(\pi_1 - 2\pi_2)} + \log^+(2(\pi_1 - 2\pi_2)) \right), & \\ & \textit{if } i\pi_i \geq (i+1)\pi_{i+1} \, \forall i \geq 1, \\ \frac{1}{(1-2\theta)\lambda}, & \textit{if } \theta < \frac{1}{2}, \end{cases}$$

*where $\lambda = \sum_{i=1}^\infty i\pi_i$ and $\theta = \frac{1}{\lambda} \sum_{i=2}^\infty i(i-1)\pi_i$.*

PROOF. The fact that $\mathscr{L}(\Psi((0,t] \times S)) = \mathscr{L}(\widetilde{\Psi}((0,t] \times S))$ and the triangle inequality imply that

$$d_{\mathrm{TV}}(\mathscr{L}(\Psi((0,t] \times S)), \mathrm{POIS}(\pi))$$

$$\leq d_{\mathrm{TV}}\bigg(\mathscr{L}(\widetilde{\Psi}((0,t]\times S)), \mathscr{L}\bigg(\int_{(0,t]\times \mathbb{Z}'_+} v\, d\xi(u,v)\bigg)\bigg)$$
$$+ d_{\mathrm{TV}}\bigg(\mathscr{L}\bigg(\int_{(0,t]\times \mathbb{Z}'_+} v\, d\xi(u,v)\bigg), \mathrm{POIS}(\pi)\bigg).$$

For the first term on the right-hand side, the basic coupling inequality and (4.1) give

$$d_{\mathrm{TV}}\bigg(\mathscr{L}(\widetilde{\Psi}((0,t]\times S)), \mathscr{L}\bigg(\int_{(0,t]\times \mathbb{Z}'_+} v\, d\xi(u,v)\bigg)\bigg) \leq 2\mathbb{P}(\widetilde{V}_1 \in S) = \frac{2\mathbb{E}(\widetilde{U}^0_{\tau^0_1-1})}{\mathbb{E}(X^0_1)}.$$

For the second term, since $\xi$ is a stationary renewal reward process, Theorem 5.1 in Erhardsson (2000b) gives a bound which equals the second term on the right-hand side in (4.3). The proof of Theorem 5.1 in Erhardsson (2000b) uses the coupling version of Stein's method for compound Poisson approximation. The last of the three bounds for the Stein constant $H_1(\pi)$ is due to Barbour and Xia (1999). □

We finally give, without proof, the lattice version of the preceding theorem.

THEOREM 4.2. *Let the conditions of Theorem* 4.1 *hold, with the following changes*: $\mathbb{R}_+$ *is replaced by* $\mathbb{Z}_+$, $\nu_\gamma$ *is the geometric distribution with mean* $\gamma^{-1}$, *and* $e^{-\gamma}$ *is replaced by* $1-\gamma$ *in the definition of* $G_s$ *for each* $s \in S$. *Then, the bound* (4.3) *remains valid, with* $\mathbb{E}((X^0_1)^2)$ *replaced by* $\mathbb{E}(X^0_1(X^0_1-1))$.

**5. Application to renewal counts.** The bound (4.3) does not at first sight seem explicit. However, by using the Markov property and solving a small number of systems of linear equations of dimension at most $N$, it is possible to express all quantities appearing in (4.3) in terms of $\gamma$, $\mu$, $\{\mathbb{E}(\zeta^S_0 I\{V^S_1 = s'\}|V^S_0 = s); (s,s') \in S\times S\}$, $\{\mathbb{E}((\zeta^S_0)^2|V^S_0 = s); s \in S\}$, $\{\int_0^\infty \sigma_s(t)\, dt; s\in S\}$ and $\{\int_0^\infty \int_u^\infty \sigma_s(t)\, dt\, du; s\in S\}$.

Below, we consider an important special case. We give a bound for the total variation distance between the distribution of the number of points in $(0,t]$ of a stationary renewal process in continuous time and a compound Poisson distribution.

By $\nu_\gamma$ we mean the exponential distribution with mean $\gamma^{-1}$.

THEOREM 5.1. *Let* $\Psi$ *be a stationary renewal point process on* $(\mathbb{R}, \mathscr{B}_\mathbb{R})$ *with generic interrenewal time* $\zeta$. *Let* $f$ *be the Radon–Nikodym derivative of the absolutely continuous part of* $\mathscr{L}(\zeta)$ *with respect to* $\nu_\gamma$. *Assume that* $\sigma : \mathbb{R}_+ \to [0,1]$ *satisfies*

$$(5.1) \qquad \sigma(t) \leq e^{-\gamma t} \inf_{x\in(t,\infty)} f(x) \qquad \forall\, t \in \mathbb{R}_+,$$

and that $G:\mathbb{R}_+ \to [0,1]$, defined by $G(t) = \sigma(t)e^{\gamma t}$, is nondecreasing and right-continuous; these conditions are satisfied if equality holds in (5.1). Let $c_0 = \gamma \int_0^\infty \sigma(t)\,dt$ and $c_1 = \gamma \int_0^\infty \int_u^\infty \sigma(t)\,dt\,du$. Assume that $c_0 > 0$. Then,

$$d_{\mathrm{TV}}(\mathscr{L}(\Psi((0,t])), \mathrm{POIS}(\pi))$$

$$\leq H_1(\pi)\frac{3t}{\mathbb{E}(\zeta)^2}\left(\frac{\mathbb{E}(\zeta) - \gamma^{-1}c_0}{c_0} + \frac{\mathbb{E}(\zeta) - c_1}{c_0} + \frac{\mathbb{E}(\zeta^2) - 2\gamma^{-1}c_1}{\mathbb{E}(\zeta)}\right.$$

$$\left. + \frac{2(\mathbb{E}(\zeta) - c_1)(\mathbb{E}(\zeta) - \gamma^{-1}c_0)}{c_0\mathbb{E}(\zeta)}\right)$$

$$+ \frac{2(\mathbb{E}(\zeta) - \gamma^{-1}c_0)}{\mathbb{E}(\zeta)},$$

where $\pi_i = \|\pi\|(1-c_0)^{i-1}c_0$ for $i \geq 1$, $\|\pi\| = tc_0/\mathbb{E}(\zeta)$, and

$$H_1(\pi) \leq \begin{cases} \left(\dfrac{1}{\|\pi\|c_0} \wedge 1\right)\exp(\|\pi\|), & \text{if } c_0 \in (0,1], \\ \dfrac{1}{\|\pi\|c_0(2c_0-1)}\left(\dfrac{1}{4\|\pi\|c_0(2c_0-1)} + \log^+(2\|\pi\|c_0(2c_0-1))\right) \wedge 1, \\ & \text{if } c_0 \in [\tfrac{1}{2},1], \\ \dfrac{c_0^2}{\|\pi\|(5c_0-4)}, & \text{if } c_0 \in (\tfrac{4}{5},1]. \end{cases}$$

PROOF. We shall compute the bound (4.3) in the case $S = \{1\}$ for a fixed $\varepsilon$, and let $\varepsilon \to 0$. All quantities appearing in (4.3) can be expressed in terms of $\gamma$, $\mathbb{E}(\zeta)$, $\mathbb{E}(\zeta^2)$, $c_0$ and $c_1$, by solving a small number of systems of linear equations. To do this, recall from Section 4 the definitions of the Markov chain $\{(\widetilde{\zeta}_i^S, \widetilde{V}_{i+1}^S); i \in \mathbb{Z}\}$ and the random sequences $\{(\widetilde{\zeta}_i^0, \widetilde{V}_{i+1}^0); i \in \mathbb{Z}\}$ and $\{(X_i^0, Y_i^0); i \in \mathbb{Z}\}$. Also, let $\tau_1 = \min\{i \in \mathbb{Z}'_+; \widetilde{V}_i^S = 0\}$ and $\tau_1^0 = \min\{i \in \mathbb{Z}'_+; \widetilde{V}_i^0 = 0\}$.

1. Clearly, $\mathbb{P}(Y_0^0 = k) = \mathbb{P}(\tau_1^0 = k+1) = \varepsilon(1-c_0)^{k-1}c_0$ for each $k \in \mathbb{Z}'_+$. In particular, $\mathbb{E}(Y_0^0) = \varepsilon/c_0$.

2. Define $h_0:\{0,1\} \to \mathbb{R}_+$ by $h_0(s) = \mathbb{E}(\sum_{i=0}^{\tau_1-1}\widetilde{\zeta}_i^S|\widetilde{V}_0^S = s)$. Conditioning on $(\widetilde{\zeta}_0^S, \widetilde{V}_1^S)$ and using the Markov property, we see that $\mathbb{E}(X_1^0) = h_0(0) = \varepsilon\gamma^{-1} + \varepsilon h_0(1)$ and $h_0(1) = \mathbb{E}(\widetilde{\zeta}_0^S|\widetilde{V}_0^S = 1) + (1-c_0)h_0(1)$. From the definition of $\widetilde{p}$ we see that $\mathbb{E}(\widetilde{\zeta}_0^S|\widetilde{V}_0^S = 1) = \mathbb{E}(\zeta) - \gamma^{-1}c_0$. It follows that $h_0(1) = (\mathbb{E}(\zeta) - \gamma^{-1}c_0)/c_0$ and $\mathbb{E}(X_1^0) = \varepsilon\mathbb{E}(\zeta)/c_0$.

3. Define $h_1:\{0,1\} \to \mathbb{R}_+$ and $h_2:\{0,1\} \to \mathbb{R}_+$ by $h_1(s) = \mathbb{E}(\sum_{i=0}^{\tau_1-1}(\widetilde{\zeta}_i^S)^2|\widetilde{V}_0^S = s)$ and $h_2(s) = \mathbb{E}(\sum_{i=0}^{\tau_1-1}\widetilde{\zeta}_i^S\sum_{j=i+1}^{\tau_1-1}\widetilde{\zeta}_j^S|\widetilde{V}_0^S = s)$, respectively. Again conditioning on $(\widetilde{\zeta}_0^S, \widetilde{V}_1^S)$ and using the Markov property, we see that

$\mathbb{E}((X_1^0)^2) = \mathbb{E}((\sum_{i=0}^{\tau_1-1} \widetilde{\zeta}_i^{\widetilde{S}})^2 | \widetilde{V}_0^{\widetilde{S}} = 0) = 2\varepsilon^2 \gamma^{-2} + \varepsilon \mathbb{E}((\sum_{i=0}^{\tau_1-1} \widetilde{\zeta}_i^{\widetilde{S}})^2 | \widetilde{V}_0^{\widetilde{S}} = 1) + 2\varepsilon^2 \gamma^{-1} h_0(1) = 2\varepsilon^2 \gamma^{-2} + \varepsilon h_1(1) + 2\varepsilon h_2(1) + 2\varepsilon^2 \gamma^{-1} h_0(1)$. Also, $h_1(1) = \mathbb{E}((\widetilde{\zeta}_0^{\widetilde{S}})^2 | \widetilde{V}_0^{\widetilde{S}} = 1) + (1-c_0)h_1(1)$, and $h_2(1) = \mathbb{E}(\widetilde{\zeta}_0^{\widetilde{S}} I\{\widetilde{V}_1^{\widetilde{S}} = 1\} | \widetilde{V}_0^{\widetilde{S}} = 1) \times h_0(1) + (1-c_0)h_2(1)$. Again, from the definition of $\widetilde{p}$ we see that $\mathbb{E}((\widetilde{\zeta}_0^{\widetilde{S}})^2 | \widetilde{V}_0^{\widetilde{S}} = 1) = \mathbb{E}(\zeta^2) - 2\gamma^{-1}c_1$, and $\mathbb{E}(\widetilde{\zeta}_0^{\widetilde{S}} I\{\widetilde{V}_1^{\widetilde{S}} = 1\} | \widetilde{V}_0^{\widetilde{S}} = 1) = \mathbb{E}(\zeta) - c_1$. It follows that $h_1(1) = (\mathbb{E}(\zeta^2) - 2\gamma^{-1}c_1)/c_0$ and $h_2(1) = (\mathbb{E}(\zeta) - c_1)(\mathbb{E}(\zeta) - \gamma^{-1}c_0)/c_0^2$. Hence,

$$\mathbb{E}((X_1^0)^2) = 2\varepsilon^2\gamma^{-2} + \frac{2(\varepsilon^2\gamma^{-1}\mathbb{E}(\zeta) - \varepsilon^2\gamma^{-2}c_0)}{c_0} + \frac{\varepsilon\mathbb{E}(\zeta^2) - 2\varepsilon\gamma^{-1}c_1}{c_0}$$
$$+ \frac{2(\mathbb{E}(\zeta) - c_1)(\varepsilon\mathbb{E}(\zeta) - \varepsilon\gamma^{-1}c_0)}{c_0^2}.$$

4. Define $h_3 : \{0,1\} \to \mathbb{R}_+$ and $h_4 : \{0,1\} \to \mathbb{R}_+$ by $h_3(s) = \mathbb{E}(\sum_{i=0}^{\tau_1-1} \sum_{j=i}^{\tau_1-1} \widetilde{\zeta}_j^{\widetilde{S}} | \widetilde{V}_0^{\widetilde{S}} = s)$ and $h_4(s) = \mathbb{E}(\sum_{i=0}^{\tau_1-1} \widetilde{\zeta}_i^{\widetilde{S}}(\tau_1 - 1 - i) | \widetilde{V}_0^{\widetilde{S}} = s)$. Yet again conditioning on $(\widetilde{\zeta}_0^{\widetilde{S}}, \widetilde{V}_1^{\widetilde{S}})$ and using the Markov property, we see that $\mathbb{E}(X_1^0 Y_0^0) = \mathbb{E}(\sum_{i=0}^{\tau_1-1} \widetilde{\zeta}_i^{\widetilde{S}}(\tau_1 - 1) | \widetilde{V}_0^{\widetilde{S}} = 0) = \varepsilon^2\gamma^{-1}\mathbb{E}(\tau_1|\widetilde{V}_0^{\widetilde{S}} = 1) + \varepsilon\mathbb{E}(\tau_1 \sum_{i=0}^{\tau_1-1} \widetilde{\zeta}_i^{\widetilde{S}} | \widetilde{V}_0^{\widetilde{S}} = 1) = \varepsilon^2\gamma^{-1}\mathbb{E}(\tau_1|\widetilde{V}_0^{\widetilde{S}} = 1) + \varepsilon h_3(1) + \varepsilon h_4(1)$. Likewise, $h_3(1) = h_0(1) + (1-c_0)h_3(1)$, and $h_4(1) = \mathbb{E}(\widetilde{\zeta}_0^{\widetilde{S}} I\{\widetilde{V}_1^{\widetilde{S}} = 1\} | \widetilde{V}_0^{\widetilde{S}} = 1)\mathbb{E}(\tau_1|\widetilde{V}_0^{\widetilde{S}} = 1) + (1-c_0)h_4(1)$. It follows that $h_3(1) = (\mathbb{E}(\zeta) - \gamma^{-1}c_0)/c_0^2$ and $h_4(1) = (\mathbb{E}(\zeta) - c_1)/c_0^2$. Hence,

$$\mathbb{E}(X_1^0 Y_0^0) = \frac{\varepsilon^2\gamma^{-1}}{c_0} + \frac{\varepsilon\mathbb{E}(\zeta) - \varepsilon\gamma^{-1}c_0}{c_0^2} + \frac{\varepsilon\mathbb{E}(\zeta) - \varepsilon c_1}{c_0^2}.$$

5. It holds that $\mathbb{E}(\widetilde{U}_{\tau_1^0-1}^0) = \mathbb{E}(\sum_{i=0}^{\tau_1-2} \widetilde{\zeta}_i^{\widetilde{S}} | \widetilde{V}_0^{\widetilde{S}} = 0) \leq \mathbb{E}(\widetilde{\zeta}_0^{\widetilde{S}} I\{\widetilde{V}_1^{\widetilde{S}} = 1\} | \widetilde{V}_0^{\widetilde{S}} = 0) + \varepsilon h_0(1) = \varepsilon^2\gamma^{-1} + (\varepsilon\mathbb{E}(\zeta) - \varepsilon\gamma^{-1}c_0)/c_0$.

We finally let $\varepsilon \to 0$ in (4.3).  $\square$

REMARK 5.1. In order to clarify what is needed to make the bound in Theorem 5.1 small, recall from Remark 3.2 the representation $\zeta = \chi(\eta_0 + \eta_1) + (1-\chi)\eta_2$, where the random variables $\chi, \eta_0, \eta_1$ and $\eta_2$ are independent, $\chi$ takes values in $\{0,1\}$ and $\eta_1$ is exponentially distributed with mean $\gamma^{-1}$. It is easy to see that $\mathbb{P}(\chi = 1) = c_0$ and $\mathbb{E}(\chi(\eta_0 + \eta_1)) = c_1$, implying that $\mathbb{E}(\zeta) - \gamma^{-1}c_0 = c_0\mathbb{E}(\eta_0) + (1-c_0)\mathbb{E}(\eta_2)$, $\mathbb{E}(\zeta) - c_1 = (1-c_0)\mathbb{E}(\eta_2)$ and $\mathbb{E}(\zeta^2) - 2\gamma^{-1}c_1 = c_0\mathbb{E}(\eta_0^2) + (1-c_0)\mathbb{E}(\eta_2^2)$.

As a consequence, assume that $c_0 \geq c > 0$ and $0 < a \leq t/\mathbb{E}(\zeta) \leq b < \infty$ (if $c > \frac{4}{5}$, the second condition is not needed). Then, the bound in Theorem 5.1 is bounded above and below by a positive constant times the expression

$$\max\left\{\frac{\mathbb{E}(\eta_0)}{\mathbb{E}(\zeta)}, \frac{\mathbb{E}(\eta_0^2)}{\mathbb{E}(\zeta)^2}, \frac{(1-c_0)\mathbb{E}(\eta_2)}{\mathbb{E}(\zeta)}, \frac{(1-c_0)\mathbb{E}(\eta_2^2)}{\mathbb{E}(\zeta)^2}\right\}.$$

REMARK 5.2. The bound given in Theorem 5.1 simplifies further if $\mathscr{L}(\zeta)$ has a Radon–Nikodym derivative $f$ with respect to $\nu_\gamma$ for some $\gamma > 0$, and $\inf_{x \in (t,\infty)} f(x) = c > 0$ for each $t \in \mathbb{R}_+$. It is then clear that we may choose $c_0 = c$ and $c_1 = \gamma^{-1} c$.

For example, assume that $\mathscr{L}(\zeta)$ is DFR (decreasing failure rate), and the failure rate has a strictly positive limit $\gamma > 0$. It then follows from Remark 4.9 in Brown (1983) that $f(x)$ decreases monotonically as $x \to \infty$ to a limit $c \geq 0$. If $c > 0$, we are in the case just described.

REMARK 5.3. Assume that $\Psi$ is a Poisson process, that is, that $\mathscr{L}(\zeta) = \nu_\gamma$ for some $\gamma > 0$. Then, from Remark 5.2, $c_0 = 1$ and $c_1 = \gamma^{-1}$, so the bound given in Theorem 5.1 is 0. The approximating distribution $\mathrm{POIS}(\pi)$ is $\mathrm{Po}(t\gamma)$.

Department of Mathematics
Royal Institute of Technology
S-100 44 Stockholm
Sweden
e-mail: torkele@math.kth.se